\documentclass{article}
\usepackage{amssymb,amsmath}
\usepackage{anysize}
\usepackage{psfrag}
\newtheorem{def.}{Definition}[section]

\newtheorem{prop}{Proposition}[section]
\newtheorem{cor}{Corollary}[section]
\newtheorem{conj}{Conjecture}[section]
\newtheorem{cl}{Claim}[section]
\newtheorem{lem}{Lemma}[section]

\numberwithin{table}{section}

\begin{document}
\title{On the minimum number of colors for knots}
       \author{Louis H. Kauffman\\
        Department of Mathematics, Statistics and Computer Science\\
        University of Illinois at Chicago\\
        851 S. Morgan St., Chicago IL 60607-7045\\
        USA\\
        \texttt{kauffman@uic.edu}\\
        and\\
        Pedro Lopes\\
        Department of Mathematics\\
        Instituto Superior T\'ecnico\\
        Technical University of Lisbon\\
        Av. Rovisco Pais\\
        1049-001 Lisbon\\
        Portugal\\
        \texttt{pelopes@math.ist.utl.pt}\\
}
\date{February 24, 2007}
\maketitle

\begin{abstract}
In this article we take up the calculation of the minimum number
of colors needed to produce a non-trivial coloring of a knot. This
is a knot invariant and we use the torus knots of type $(2, n)$ as
our case study. We calculate the minima in some cases. In other
cases we estimate upper bounds for these minima leaning on the
features of modular arithmetic. We introduce a sequence of
transformations on colored diagrams called Teneva transformations.
Each of these transformations reduces the number of colors in the
diagrams by one (up to a point). This allows us to further
decrease the upper bounds on these minima. We conjecture on the
value of these minima. We apply these transformations to rational
knots.
\end{abstract}

\bigbreak

Keywords: Knots, colorings, colors, Teneva transformations

\bigbreak

2000 MSC: 57M27

\bigbreak

\section{Introduction} \label{sect:intro}

\noindent

The colorings we are concerned with are the so-called Fox
colorings, \cite{Fox, lhKauffman}. Given a knot diagram and an integer $r$ we
consider the integers $0, 1, 2, \dots , r-1 $ mod $r$, whose set
will be denoted $\mathbb{Z}\sb{r}$. We assign one integer (call it
a color) to each arc of the diagram so that at each crossing the
sum of the integers at the under-arcs minus twice the integer at
the over-arc equals zero mod $r$. In this way, we set up a system
of equations over $\mathbb{Z}\sb{r}$. Each of the solutions of
this system of equations is called an $r$-coloring of the knot
under consideration. There are always at least $r$ solutions,
regardless of the knot diagram or the integer $r$ we are choosing;
each of these $r$ solutions is obtained by assigning the same
color to each arc of the diagram. These are called the trivial
solutions. A one-to-one correspondence of the solutions of the
systems of equations before and after the performance of each of
the Reidemeister moves is presented in \cite{pLopes}. In this
correspondence, the trivial solutions go over to trivial solutions
and non-trivial solutions go over to non-trivial solutions. In
particular, the number of solutions for each $r$ (also referred to
as the number of $r$-colorings), is an invariant of the knot under
study. Given an integer $r>1$ and a link $K$, we let $\# col\sb{r}
K$ stand for the number of $r$-colorings of $K$.

The efficiency of this invariant in distinguishing prime knots up
to ten crossings is illustrated in \cite{DL}. The invariant we
there associate to each knot is in fact the {\it color spectrum}
of the knot i.e., the sequence of numbers of $r$-colorings. We
remark that the notion of {\it quandles} generalizes the notion of
$r$-colorings, see \cite{dJoyce, {sMatveev}}.

The importance of knowing the color spectrum of the knot $K$ under
consideration, is that it tells us immediately whether $K$ is
interesting or not for our current research and whether there is
significant topological information in the spectrum. In fact, for
any $r$, the number of trivial colorings is $r$. Hence if the
number of colorings for a given $r$ is greater than $r$, there
exist then non-trivial colorings for the knot under consideration.
As an example, the trefoil exhibits nine $3$-colorings. It has
then non-trivial $3$-colorings.

\bigbreak

\begin{def.}[Minimum number of colors]
Given an integer $r>1$, assume there are non-trivial $r$-colorings
on a given knot $K$. Assume further that $D$ is a diagram of $K$.
We let $n\sb{r, K}(D)$ denote the minimum number of distinct
colors assigned to the arcs of $D$ it takes to produce a
non-trivial $r$-coloring on $D$. We denote
$\mathrm{mincol}\sb{r}K$ the minimum of these minima over all
diagrams $D$ of $K$:
\[
\mathrm{mincol}\sb{r}K := \min \{ n\sb{r, K}(D) \; | \; D
\textnormal{ is a diagram of } K \}
\]
For each $K$, we call $\mathrm{mincol}\sb{r}K$ the {\bf minimum
number of colors of K, mod $r$}. In the sequel, we will drop the
``mod $r$'' whenever it is clear which $r$ we are referring to.
Note that $\mathrm{mincol}\sb{r}K$ is tautologically a topological
invariant of $K$.
\end{def.}

Apparently, in order to calculate $\mathrm{mincol}\sb{r}K$, we
have to consider a diagram of $K$, and find the minimum number of
colors it takes to construct a non-trivial $r$-coloring. This
operation should then be repeated for all diagrams of $K$,
presenting the minimum for each of these diagrams. Finally, the
minimum of these minima is the $\mathrm{mincol}\sb{r}(K)$. In this
article we present techniques that allow us to calculate
$\mathrm{mincol}\sb{r}(K)$ in infinitely many cases and in other
cases to estimate its upper bounds. We regard the torus knots of
type $(2, n)$ as our case study. The features of modular
arithmetic allow us to calculate $\mathrm{mincol}\sb{r}T(2, n)$
exactly for some combinations of $n$ and $r$. The introduction of
certain transformations on diagrams (Teneva transformations)
allows us to better estimate the upper bounds on
$\mathrm{mincol}\sb{r}T(2, n)$ for other combinations of $n$ and
$r$.

\bigbreak

This invariant was first introduced in \cite{hk}. Considerations
of the authors about it led them to set forth the Kauffman-Harary
Conjecture which has already been proven to be true for rational
knots, \cite{kl-t}, and for Montesinos links, \cite{aps}.

\begin{conj} [Kauffman-Harary, \cite{hk}]\label{conj:KH} Let
$p$ be a prime integer and assume $K$ is an alternating knot of
determinant $p$. Then, any non-trivial $p$-coloring on any minimal
diagram of $K$ assigns different colors to different arcs of the
diagram.
\end{conj}

We note that the Kauffman-Harary Conjecture deals only with a
specific $r$ per knot. Moreover, the knots under consideration in
this conjecture are all alternating knots of prime determinant
(and the specific $r$ is precisely this determinant, for each of
these knots). In this article we adopt a broader point of view by
not specifying $r$. In fact, we would like to develop
computational tools that would allow us to calculate the
$\mathrm{mincol}\sb{r}(K)$ for any $r$, and for any $K$. Here we
content ourselves on studying the class of torus knots of type
$(2, n)$.

\bigbreak

Given two positive integers $l$ and $m$ we let $(l, m)$ stand for
the greatest common divisor of $l$ and $m$ and $\langle l, m
\rangle$ stand for $1$ if $(l, m)=1$, and for the least common
prime divisor of $l$ and $m$ otherwise.

\bigbreak

In this article we prove the following Theorem:

\bigbreak

{\bf Main Theorem} \quad Suppose $r$ and $n$ are positive integers
such that $(n, r)>1$.
\begin{itemize}
\item If $\langle n, r \rangle \in \{ 2, 3 \}$ then
\[
\mathrm{mincol}\sb{r}T(2, n)=\langle n, r \rangle
\]
\item If $\langle n, r \rangle = 5$ then
\[
\mathrm{mincol}\sb{r}T(2, n)=4
\]
\item If $\langle n, r \rangle = 2k+1$ (for some integer $k>2$)
then
\[
3 < \mathrm{mincol}\sb{r}T(2, n) \leq k+2
\]
\end{itemize}

\bigbreak

Moreover, we conjecture that

\bigbreak

{\bf Conjecture}  \quad Suppose $r$ and $n$ are positive integers
such that $(n, r)>1$. Let $\langle n, r \rangle =2k+1$, for some
integer $k>2$. Then:
\[
\mathrm{mincol}\sb{r}T(2, n) = k+2
\]

\bigbreak

\bigbreak

The organization of this article is as follows. In Section
\ref{sect:Spectral} we calculate the number of $r$ colorings for
each $T(2, n)$ and estimate upper bounds on the minimum number of
colors for each $T(2, n)$. In Section \ref{sect:proof} we prove
the Main Theorem. In Section \ref{sect:applications} we present
applications to rational knots. In Section \ref{sect:final} we
consider directions for further work.

\section{The number of $r$-colorings of $T(2, n)$. \\
Upper bounds on minimum numbers of colors.}\label{sect:Spectral}

\noindent

\begin{figure}[h!]
    \psfrag{a1}{\huge $a$}
    \psfrag{b1}{\huge $b$}
    \psfrag{a2}{\huge $2b-a$}
    \psfrag{a3}{\huge $3b-2a$}
    \psfrag{a4}{\huge $4b-3a$}
    \centerline{\scalebox{.50}{\includegraphics{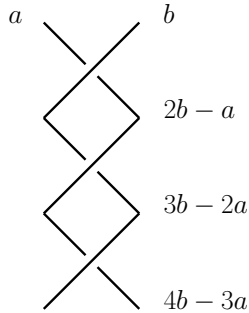}}}
    \caption{Coloring assignment to $\sigma\sb 1\sp 3$}\label{Fi:sigma1n3kh}
\end{figure}

\begin{figure}[h!]
    \psfrag{a1}{\huge $a$}
    \psfrag{b1}{\huge $b$}
    \psfrag{a2}{\huge $2b-a$}
    \psfrag{d}{\huge $\; \: \vdots$}
    \psfrag{an}{\huge $nb-(n-1)a$}
    \psfrag{an+1}{\huge $(n+1)b-na$}
    \psfrag{an+2}{\huge $(n+2)b-(n+1)a$}
    \centerline{\scalebox{.50}{\includegraphics{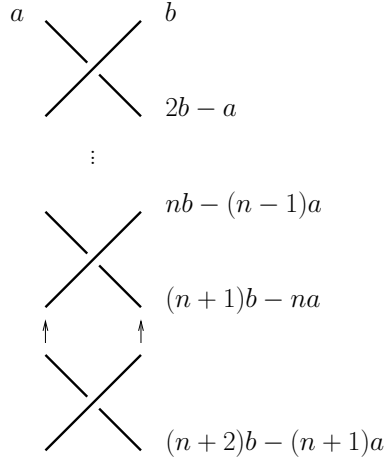}}}
    \caption{Juxtaposing the $(n+1)$-th crossing}\label{Fi:sigma1nnkh}
\end{figure}

\begin{prop}\label{prop : spectral}For each integer $n\geq2$,
\[
\# col\sb{r} T(2, n) = (n, r)r
\]
for all $r\geq 2$.\end{prop}Proof: Consider the diagram of $T(2,
n)$ given by the closure of $\sigma\sb {1}\sp {n}$ ($\sigma\sb {1}
\in B\sb {2}$) - see Figure \ref{Fi:sigma1n3kh} for the $n=3$
instance and \cite{Birman} for further information on the
$\sigma\sb{i}$'s.
Assume further that
$\widehat{\sigma\sb{1}\sp {n}}$ is endowed with an $r$-coloring
for some integer $r\geq 2$. In order to set up the system of
equations whose solutions are the $r$-colorings of $T(2, n)$, we
will state and prove the following claim:
\begin{cl} \label{cl:colonbraid} Let $n$ be an integer greater than
$2$ and assume $\sigma\sb{1}\sp {n}$ is
endowed with an $r$-coloring for some integer $r\geq 2$. If the
top segments of the braid $\sigma\sb {1}\sp {n}$ $(\sigma\sb {1}
\in B\sb {2})$ are colored $($from left to right$)$ with $a$,
$b\in R\sb r$ then, for each $i=1,...,n$, the segment emerging
from the $i$-th crossing will bear color $(i+1)b-ia$ $($counting
crossings from top to bottom$)$.\end{cl}Proof: By induction on
$n$. The $n=3$ instance is clear by inspection of figure
\ref{Fi:sigma1n3kh}. Assume claim is true for a specific $n\in
\mathbb{N}$. Juxtaposing another crossing to the previous $n$ (see
figure \ref{Fi:sigma1nnkh}), it is easy to see that the induction
step follows. $\hfill \blacksquare $

We now return to the proof of the Proposition. For $n$ crossings,
the labellings on the bottom segments become $nb-(n-1)a$ and
$(n+1)b-na$ (from left to right). In order for these $a$ and $b$
to stand for an $r$-coloring, the bottom left labelling has to
equal the top left labelling and the bottom right labelling has to
equal the top right labelling i.e., $nb-(n-1)a=\sb{r} a$ and
$(n+1)b-na=\sb{r} b$ which simplify both to
\[
n(b-a)=\sb{r} 0
\]
So, $b-a$ has to provide the largest factor of $r$ which is
relatively prime to $n$ for the equation to hold; this factor is
$\frac{r}{(n, r)}$. Thus $b-a$ equals $\frac{r}{(n, r)}$ or one of
its multiples in $\{ 1, 2, \dots , r  \}$ i.e., $b-a$ can be any
of the following:
\[
\frac{r}{(n, r)}, \quad 2\frac{r}{(n, r)}, \quad 3\frac{r}{(n,
r)}, \quad \dots \quad , \quad  \bigl( (n, r) - 1 \bigr)
\frac{r}{(n, r)}, \quad (n, r)\frac{r}{(n, r)}
\]
Finally,  there are $r$ ordered pairs $(a, b)$ compatible with
each of the above $(n, r)$ possibilities. This statement is
justified by the following Claim.
\begin{cl} For any $i\in \{ 0, 1, 2, \dots , r-1  \}$, there are
exactly $r$ pairs $(a, b)$ which are solutions to $b-a=\sb{r} i$
$(a, b\in \{ 0, 1, 2, \dots , r-1 \})$.\end{cl}Proof: For each assignment of a value to $a$, $b$ is uniquely specified by the formula $b=\sb{r} a + i$. Since $r$ values can be assigned to $a$, the proof is complete. $\hfill
\blacksquare $

Resuming the proof of the Proposition, there are $(n, r)$
possibilities for the $b-a$ so that the equation $n(b-a)=\sb{r} 0$
holds. Each of these possibilities can be realized in $r$ distinct
ways. There are then $(n, r)r$ solutions i.e., $(n, r)r$
$r$-colorings of $T(2, n)$. $\hfill \blacksquare $

\bigbreak

\begin{cor}\label{cor : nontriv}  There are nontrivial colorings
if, and only if, $(n, r)\neq 1$. \end{cor} Proof: Omitted. $\hfill
\blacksquare $

\bigbreak

We recall that, given two positive integers $l$ and $m$, $\langle
l, m \rangle $ stands for $1$ if $(l, m)=1$ and stands for the
least common prime factor of $l$ and $m$, otherwise.

\begin{prop}\label{prop : mincor} Let $n$ and $r$ be positive integers
greater than $1$. If $(n, r)=1$ there are only trivial
$r$-colorings of $T(2, n)$. If $(n, r)>1$,
\[
\mathrm{mincol}\sb{r}T(2, n) \leq \langle n, r \rangle
\]
\end{prop}Proof: That for $( n, r ) = 1$ there are only trivial
colorings is just a rephrasing of Corollary \ref{cor : nontriv}.

If $( n, r ) > 1$ set
\[
p = \langle n, r \rangle
\]
The set
\[
R\sb{r}\sp{p}:=\left\{  0, \; \frac{r}{p} , \; 2\: \frac{r}{p} ,\:
\dots \: , \; (p-1)\: \frac{r}{p} \right\}
\]
endowed with the $a\ast b := 2b-a$ (mod $r$) operation
(\cite{dJoyce, sMatveev}) is algebraically closed. As a matter of
fact, it is an algebraically closed substructure of $\{ 0, 1, 2,
\dots r-1 \}$ endowed with the same operation. It is equivalent to
the set $\{ 0, 1, \dots p-1  \}$ endowed with the operation $a\ast
b := 2b-a$ (mod $p$).

\begin{figure}[h!]
    \psfrag{a}{\huge $a$}
    \psfrag{b}{\huge $b$}
    \psfrag{2b-a}{\huge $2b-a$}
    \psfrag{3b-2a}{\huge $3b-2a$}
    \psfrag{A}{\huge $pb-(p-1)a=a$}
    \psfrag{B}{\huge $(p+1)b-pa=b$}
    \psfrag{...}{\huge $\vdots $}
    \psfrag{p xings}{\huge $\sigma\sb 1\sp p$ }
    \psfrag{ith xing}{\huge $i$-th xing}
    \psfrag{i+p-1th xing}{\huge $(i+p-1)$-th xing}
    \centerline{\scalebox{.50}{\includegraphics{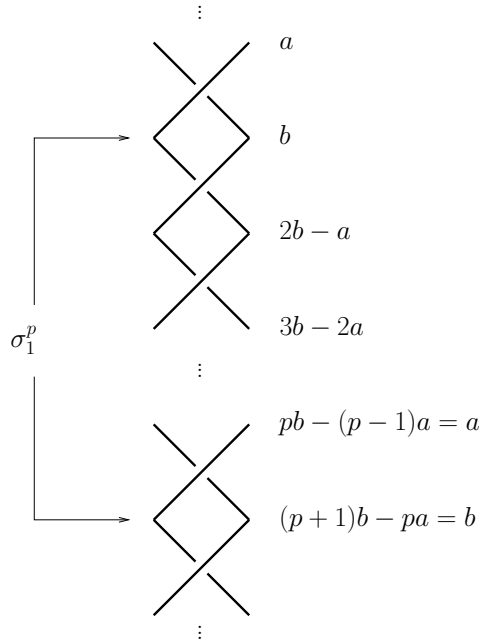}}}
    \caption{Coloring of $\sigma\sb{1}\sp{p}$ inside coloring of
    $\sigma\sb{1}\sp{n}$}\label{Fi:sigma1pnp}
\end{figure}

Also note that, in this $p=\langle n, r \rangle > 1$ case, the
$\sigma\sb{1}\sp{n}$ braid ($\sigma\sb{1}\in B\sb{2}$), whose
closure gives a knot diagram for $T(2, n)$, can be regarded as a
product of $\frac{n}{p}$ $\sigma\sb{1}\sp{p}$'s:
\[
\underbrace{\sigma\sb 1\sp p \cdot \dots \cdot \sigma\sb 1\sp
p}_{\text{ $\frac{n}{p}$ factors }}
\]
This reflects on the coloring equation:
\[
0 =\sb{r} n(b-a) =\sb{r} \frac{n}{p} \cdot p(b-a)
\]
Hence, any two $a$, $b$ from $R\sb{r}\sp{p}$ yield an $r$-coloring
of $\widehat{\sigma\sb 1\sp n}$. We have, thus, $p\sp{2}$
$r$-colorings with colors from $R\sb{r}\sp{p}$. Note, also, that
any of these $a$, $b$ from $R\sb{r}\sp{p}$ yield an $r$-coloring
of each $\widehat{\sigma\sb 1\sp p}$, since they satisfy
\[
0 =\sb{r}  p(b-a)
\]
We can, thus, regard some non-trivial $r$-colorings of
$\widehat{\sigma\sb 1\sp n}$ as stackings of $\frac{n}{p}$
$p$-colorings of $\widehat{\sigma\sb 1\sp p}$ (see Figure
\ref{Fi:sigma1pnp}).
Now, each of these non-trivial $p$-colorings of
$\widehat{\sigma\sb 1\sp p}$ uses exactly $p$ distinct colors from
$R\sb{r}\sp{p}$ thanks to the following claim:

\begin{cl}\label{cl : prime} If $p$ is prime then each non-trivial
$p$-coloring of $\widehat{\sigma\sb 1\sp p}$ uses exactly $p$
distinct colors.\end{cl}Proof: Pick $a$, $b\in \{ 0, 1, \dots ,
p-1 \}$, distinct. Suppose there exist $i, j\in \{ 0, \dots , p-1
\}$ such that $ib-(i-1)a=\sb{p} jb-(j-1)a$. Then
$(i-j)(b-a)=\sb{p} 0$ which is equivalent to saying that $i=j$,
for $\mathbb{Z}\sb p$ is a field.$\hfill \blacksquare $

Hence, given positive integers $n$ and $r$ such that $(n, r)>1$ we
constructed a non-trivial $r$-coloring of $T(2, n)$ which uses
exactly $\langle n, r \rangle $ colorings. $\langle n, r \rangle $
is then an upper bound for $\mathrm{mincol}\sb{r}T(2, n)$. This
concludes the proof. $\hfill \blacksquare $

\section{Proof of the Main Theorem.}\label{sect:proof}

\noindent

In this Section we prove the Main Theorem. In Subsection
\ref{subsect:23} we prove the $\langle n, r \rangle \in \{ 2, 3
\}$ instance. In Subsection \ref{subsect:TT} we prove $k+2$ is an
upper bound on the minimum number of colors of the $\langle n, r
\rangle = 2k+1$ instance after introducing the Teneva
transformations. In Subsection \ref{subsect:5} we prove the
$\langle n, r \rangle = 5$ instance and that $3$ is a lower bound
on the minimum number of colors of the $\langle n, r \rangle =
2k+1$ instance.

\subsection{The $\langle n, r \rangle \in \{ 2, 3 \}$ instance.}\label{subsect:23}

\begin{prop} For even positive integers $n$ and $r$
\[
\mathrm{mincol}\sb{r}\: T\bigl( 2, n \bigr) = 2
\]
\end{prop} Proof: For even positive integers $n$ and $r$, $ \langle n, r \rangle = 2$. Then, by
Proposition \ref{prop : mincor},
\[
\mathrm{mincol}\sb{r}\: T\bigl( 2, n \bigr) \leq 2
\]
Since a nontrivial coloring has to use at least two distinct
colors, the result follows. $\hfill \blacksquare $

\bigbreak

We remind again the reader that we call trivial coloring any
coloring which assigns the same color to each arc of the diagram
under study. We remark that this implies that a trivial knot of
more than one component can be assigned non-trivial colorings. In
order to see this, consider a diagram of this knot where one
component lies in a neighborhood which is disjoint from a
neighborhood which contains the rest of the diagram. Color the
singled out component with color $a$ and the rest of diagram with
color $b (\neq a)$ - this is a non-trivial coloring of this
trivial knot. Note further that this type of phenomenon does not
occur for knots. In fact, a trivial knot has a diagram with no
crossings. Any coloring of this diagram can only have one color -
hence any other of this knot's diagrams is colored with only one
color.

\begin{prop}\label{prop : mincol3} For positive integers $n$ and $r$
such that $ \langle n, r \rangle = 3$,
\[
\mathrm{mincol}\sb{r}\: T\bigl( 2, n \bigr) = 3
\]
\end{prop} Proof: Since $\langle n, r \rangle =3$, then there is an
$r$-coloring of $T(2, n)$ with as few as three colors by the proof
of Proposition \ref{prop : mincor}. We next prove that two
distinct colors are not enough to produce a non-trivial coloring
in the $ \langle n, r \rangle = 3$ case. We consider two
possibilities: odd $r$ and even $r$.

Suppose $r$ is odd and $a, b\in \{ 0, 1, \dots , r-1  \}$. If
$b=2b-a$ then $a=b$; if $a=2b-a$ then $2b=2a$ which is equivalent
to saying $a=b$ since $r$ is odd (hence $2$ is invertible). Thus
if we choose distinct $a, b$ from $\{ 0, 1, \dots , r-1  \}$ then
$\# \{ a, b, 2b-a \} = 3$ and so any non-trivial $r$-coloring of
$T(2, n)$, for odd $r$, has at least three colors. This concludes
the proof for odd $r$.

Suppose $r$ is even; then $n$ is odd for otherwise $ \langle n, r
\rangle = 2$. Assume there is an $r$-coloring of a diagram $D$ of
$T(2, n)$ which uses exactly two distinct colors, say $a, b$. At
some crossing of the diagram the two colors meet. The
possibilities for the color on the emerging arc are as follows.
Either $b=2b-a$ i.e., $a=b$ which is contrary to the assumption;
or $a=2b-a$ i.e., $a=b+\frac{r}{2}$. Hence,
$2a-b=2(b+\frac{r}{2})-b=b$ mod $r$. Thus, for even $r$, $b$ and
$a=b+\frac{r}{2}$ generate an algebraically closed structure with
respect to the $x\ast y :=2y-x$ operation (\cite{dJoyce,
sMatveev}). This structure is formed precisely by $b$ and
$a=b+\frac{r}{2}$. Consider again the diagram $D$ endowed with the
indicated coloring which uses only the two colors $a$ and $b$. By
performing Reidemeister moves and consistently changing the colors
after each move (cf. \cite{pLopes}) we obtain a new diagram of
$T(2, n)$ endowed with a coloring which uses only the two colors
$a$ and $b$ - because $\{ a, b \}$ constitute an algebraically
closed set with respect to the  $\ast $-operation, as was seen
above.

In particular, we could transform $D$ into $\widehat{\sigma\sb
{1}\sp {n}}$ upon performance of Reidemeister moves and the
associated coloring (obtained by consistently changing the
colorings after each Reidemeister move) would use exactly two
colors. We now prove that $\widehat{\sigma\sb {1}\sp {n}}$ cannot
be colored with only two colors - because $n$ is odd. As a matter
of fact, starting with distinct $a$ and $b$ from $\{ 0, 1, \dots ,
r-1 \}$ at the top segments of $\sigma\sb {1}\sp {n}$ (from left
to right) we obtain, using induction, $b, a$ (from left to right)
after an odd number of crossings. This concludes the proof.
$\hfill \blacksquare $

\bigbreak

We remark that the proof of Proposition \ref{prop : mincol3}
yields results stronger than the statement of the Proposition.
Specifically,

\begin{cor}\label{cor:mincolcor1} Let $r$ be an integer greater
than $1$ and assume $K$ is not splittable.
\begin{itemize}
\item If $r$ is odd then $\mathrm{mincol}\sb{r}K > 2$ \item If $r$
is even then,
\begin{itemize}
\item a specific diagram of $K$ admits an $r$-coloring with
exactly two colors if, and only if, any other diagram of $K$
admits an $r$-coloring with exactly two colors
\end{itemize}
\end{itemize}
\end{cor}Proof: Omitted. $\hfill \blacksquare $

\bigbreak

In particular,
\begin{cor}\label{cor:mincolcor2} If $\langle n, r \rangle $ is
an odd prime then:
\[
\mathrm{mincol}\sb{r}T(2, n) > 2
\]
\end{cor}Proof: If $r$ is odd, the statement of this Corollary is a particular case of the
first statement of Corollary \ref{cor:mincolcor1}. If $r$ is even,
then $n$ is odd. In particular, $T(2, n)$ is a knot and  the
second statement of Corollary \ref{cor:mincolcor1} applies.
Repeating an argument used in the proof of Proposition \ref{prop :
mincol3} we see that the closure of $\sigma\sb{1}\sp{n}$ cannot
have an $r$-coloring with just two colors, since $n$ is odd.
Whence no diagram of $T(2, n)$ can have an $r$-coloring with just
two colors. This concludes the proof. $\hfill \blacksquare $

\subsection{Teneva transformations.}\label{subsect:TT}

\noindent

The results we obtained on $\mathrm{mincol}\sb{r}T(2, n)$ so far,
relied on the features of modular arithmetic. We will now come up
with better estimates for $\mathrm{mincol}\sb{r}T(2, n)$ by making
use  of Reidemeister moves. In particular, we will obtain diagrams
endowed with non-trivial colorings that use less colors than the
ones considered so far, although these diagrams have more arcs than
the $\widehat{\sigma\sb{1}\sp{n}}$'s. In order to obtain these
diagrams we will use what we call Teneva transformations. This is
a particular sequence of Reidemeister moves, starting from the
$\widehat{\sigma\sb{1}\sp{n}}$ diagram of $T(2, n)$ endowed with a
non-trivial coloring and consistently coloring the diagrams after
each move, in the sense introduced in \cite{pLopes}. This
formalizes and generalizes a particular case due to Irina Teneva
presented in \cite{hk}.

\bigbreak

Specifically, we now establish that, for any odd prime $p=2k+1$
(for some positive integer $k$),
\[
\mathrm{mincol}\sb{p}T(2, p) \leq k+2
\]
In order to do this we will prove that, for any positive integer
$k$, $\mathrm{mincol}\sb{2k+1} T(2, 2k+1) \leq k+2$. Of course,
for non-prime $2k+1$, Proposition \ref{prop : mincor} presents a
strictly smaller upper bound but in this way we are able to use
induction on $k$ thus establishing the $k+2$ upper bound also for
prime $p=2k+1$. This is a better estimate than the Proposition's,
in the prime $p$ situation. We remark, in passing, that, for any
positive integer $k$, $k+2$ is also an upper bound for
$\mathrm{mincol}\sb{2k}T(2, 2k)$; again in this case, Proposition
\ref{prop : mincor} presents a strictly smaller upper bound.

In order to prove that $k+2$ is an upper bound, as referred to
above, we  consider the diagram of the torus knot $T(2, 2k+1)$
(for $k>1$) as given by the braid closure of $\sigma\sb {1}\sp
{2k+1}$ ($\sigma\sb {1} \in B\sb {2}$) and endowed with a
non-trivial $(2k+1)$-coloring. This coloring is represented by
assigning $a$, $b\in \{  0, 1, \dots , 2k \}$ ($a \neq b$) to the
top segments of the braid $\sigma\sb {1}\sp {2k+1}$ in the usual
manner. Since $2k+1$ is the number of arcs of $\widehat{\sigma\sb
{1}\sp {2k+1}}$ there are, at most, $2k+1$ colors in this
coloring. If $a=0$ and $b=1$ then there are exactly $2k+1$ colors.
If $2k+1$ is prime then there are exactly $2k+1$ colors whenever
$a\neq b$, thanks to Claim \ref{cl : prime}. We, thus, assume
that, for each $k$, the coloring in point uses $2k+1$ colors. We
then perform $1+k$ Reidemeister moves on the diagrams,
consistently changing the coloring assignments after each move
(cf. \cite{pLopes}), eventually obtaining a coloring assignment
with $k+2$ colors.

\begin{figure}[h!]
    \psfrag{a1}{\huge $a$}
    \psfrag{b1}{\huge $b$}
    \psfrag{a2}{\huge $2b-a$}
    \psfrag{a3}{\huge $3b-2a$}
    \psfrag{a4}{\huge $4b-3a$}
    \psfrag{1}{\huge $I$}
    \psfrag{3}{\huge $III$}
    \centerline{\scalebox{.50}{\includegraphics{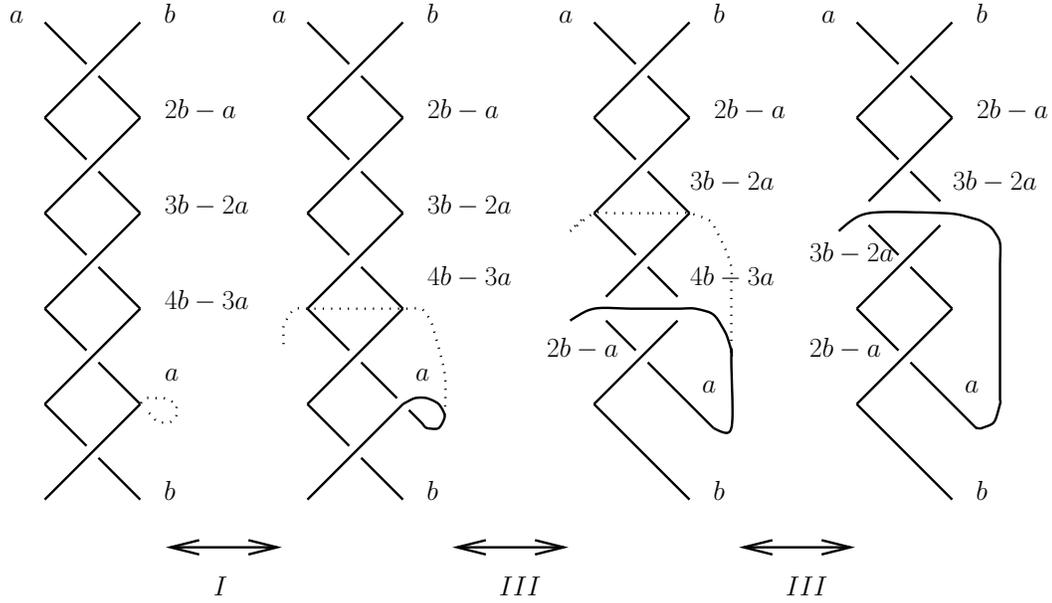}}}
    \caption{Colored T(2, 5): from 5 to 4 colors}\label{Fi:t2n5kh}
\end{figure}

We exemplify for $k=2$ starting from the diagram of $T(2, 5)$
given by the closure of $\sigma\sb {1}\sp {2\cdot 2+1}$ and
endowed with an $r$-coloring as shown (see Figure
\ref{Fi:t2n5kh}). The dotted lines indicate where the arc in point
is taken to by the next move. There are $1+2=3$ moves. The first
one is a type I Reidemeister move on the arc labeled $a$ at the
bottom right of the first diagram in figure \ref{Fi:t2n5kh}; the
two other moves are type III Reidemeister moves. The first of
these two type III Reidemeister moves moves the arc labeled $a$
(which stems from the bottom left of the second diagram) over a
crossing above it in the manner indicated in figure
\ref{Fi:t2n5kh} (dotted lines in the second diagram). Color $2b-a$
is introduced with this move but since it was already part of the
coloring, the number of colors remains the same. The second type
III Reidemeister move moves the arc we just referred to over the
crossing right above it. With this move color $3b-2a$ is
introduced and color $4b-3a$ is removed. Since color $3b-2a$ was
already part of the coloring and there was only one arc assigned
color $4b-3a$, this move reduces the number of colors by one. We
end up thus with $5-1=k+2$ colors for $k=2$ as announced.

For a general integer $k>1$, there are $1+k$ moves on $\sigma\sb
{1}\sp {2k+1}$ with the first three as indicated for $T(2, 5)$, in
figure \ref{Fi:t2n5kh}. As of the third, each move will be
analogous to the previous one, the arc in point being pulled up
over one crossing at a time. As with $T(2, 5)$, the first two
moves will not affect the number of colors used. The third move
will remove color $(2k)b-(2k-1)a$, the fourth move will remove
color $(2k-1)b-(2k-2)a$, ... , the $k$-th move will remove color
$(k+3)b-(k+2)a$, and the $(k+1)$-th move will remove color
$(k+2)b-(k+1)a$.  Analogous remarks apply for the even $n=2k$
case. We then state:

\begin{prop}[Teneva reduction]\label{prop : 2k+1} For any integer $k>1$,
\[
\mathrm{mincol}\sb{2k+1}T(2, 2k+1) \leq k+2
\]
\end{prop}

Before embarking on the proof of this proposition we will first
state and prove a technical lemma which will help us deal with it.

\begin{lem}[Teneva transformation]\label{lem : remint} For each
integer $n\geq 3$, consider $\sigma\sb {1}\sp {n}$ $(\sigma\sb {1}
\in B\sb {2})$ endowed with an $r$-coloring as in Claim \ref{cl:colonbraid}. Then, $n-1$ type III Reidemeister moves are performed, taking the bottom left arc
of the braid over each crossing above it until the top of the
braid is reached. These type III Reidemeister moves are preceded
by a type I Reidemeister move which prepares the setting for the
subsequent moves $($see left-most braid in figure \ref{Fi:t2n5kh},
for the $n=5$ instance$)$. The diagrams are consistently colored
after each move.

The first type III Reidemeister move introduces color
$(n+2)b-(n+1)a$ on the left-hand side of the diagram. For $2\leq i
\leq n-1$, the $i$-th type III Reidemeister move removes color
$(n+1-i)b-(n-i)a$ from the right-hand side of the diagram and
introduces color $(n+1+i)b-(n+i)a$ on the left-hand side of the
diagram.\end{lem}

\begin{figure}[h!]
    \psfrag{x0}{\huge $a$}
    \psfrag{x1}{\huge $b$}
    \psfrag{x2}{\huge $2b-a$}
    \psfrag{x3}{\huge $3b-2a$}
    \psfrag{y0}{\huge $3b-2a$}
    \psfrag{y1}{\huge $4b-3a$}
    \psfrag{y2}{\huge $5b-4a$}
    \psfrag{y3}{\huge $6b-5a$}
    \psfrag{1}{\huge $I$}
    \psfrag{3}{\huge $III$}
    \centerline{\scalebox{.50}{\includegraphics{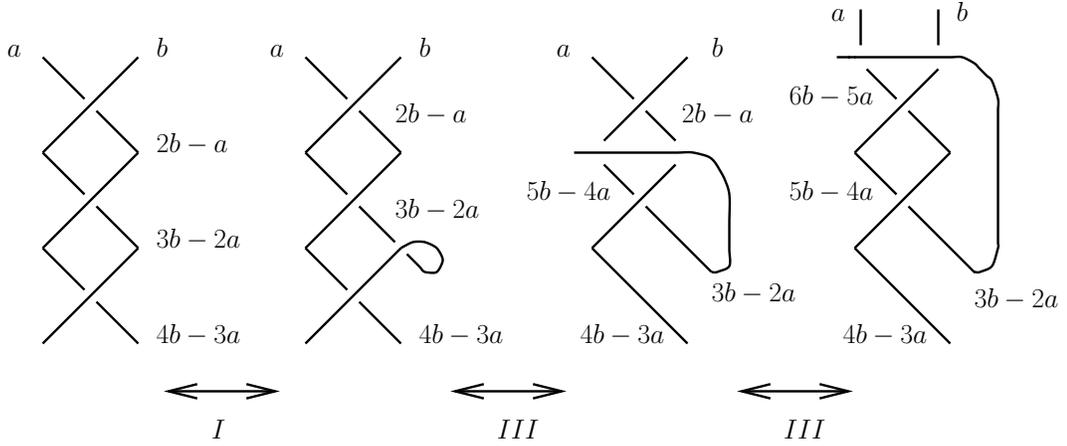}}}
    \caption{The $n=3$ instance}\label{Fi:t2n3kh}
\end{figure}

Proof: We use induction on $n$. For $n=3$, see figure
\ref{Fi:t2n3kh}.

\begin{figure}[h!]
    \psfrag{x0}{\huge $a$}
    \psfrag{x1}{\huge $b:=a'$}
    \psfrag{x2}{\huge $2b-a:=b'$}
    \psfrag{x3}{\huge $3b-2a = 2b'-a'$}
    \psfrag{d}{\huge $\; \: \vdots$}
    \psfrag{xn-1}{\huge $nb-(n-1)a=(n-1)b'-(n-2)a'$}
    \psfrag{xn}{\huge $(n+1)b-na=nb'-(n-1)a'$}
%   \psfrag{y0}{\huge $y\sb 0$}
    \psfrag{y1}{\huge $(n+2)b-(n+1)a=(n+1)b'-na'$}
    \centerline{\scalebox{.50}{\includegraphics{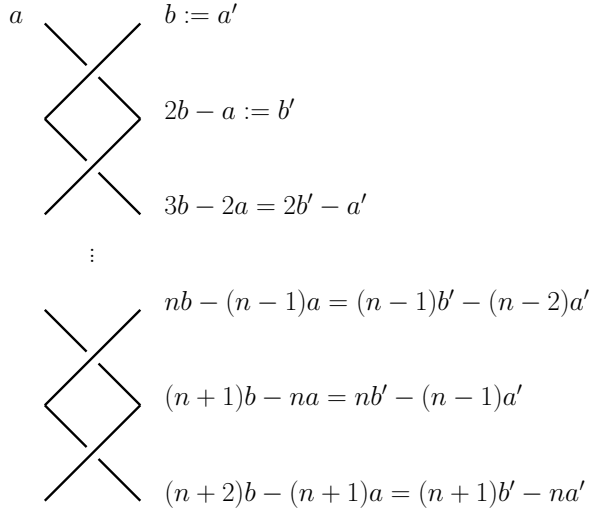}}}
    \caption{The assignment of colors to the $\sigma\sb{1}\sp{n+1}$
    braid}\label{Fi:sigma1n(n+1)lemmakh}
\end{figure}

Now assume the claim is true for a specific integer $n\geq 3$ and
consider the $n+1$ situation (see figure
\ref{Fi:sigma1n(n+1)lemmakh}).
Neglecting the top $\sigma \sb 1$, there is a string of $n$
$\sigma \sb 1$'s with the top segments colored with $a'=b$ and $b'=2b-a$. Hence we can apply the induction hypothesis to say that the first type III Redemeister
move introduces color $(n+2)b'-(n+1)a' \; \big( =
(n+3)b-(n+2)a\big) $ and, for $2\leq i \leq n-1$, the $i$-th type
III Reidemeister move removes color $(n+1-i)b'-(n-i)a' \;  \big( =
(n+2-i)b-(n+1-i)a\big) $ from the right-hand side of the diagram
and introduces color $(n+1+i)b'-(n+i)a' \; \big( =
(n+2+i)b-(n+1+i)a\big) $. The $n$-th type III Reidemeister move
removes color $2b-a$ and introduces color $(2n+2)b-(2n+1)a$.
Hence, the claim follows.$\hfill \blacksquare $

Proof (of Proposition \ref{prop : 2k+1}): For some integer $k>1$,
consider a non-trivial $(2k+1)$-coloring of
$\widehat{\sigma\sb{1}\sp{2k+1}}$, using $2k+1$ distinct colors. Lemma \ref{lem : remint}, specialized to $n=2k+1$ and reading colors mod $2k+1$, implies
that, for $2 \leq i \leq 2k$, the $i$-th type III Reidemeister
move removes color $(2k+2-i)b-(2k+1-i)a$ and introduces color
$(i+1)b-ia$. Write the set of colors in the original coloring of
$\sigma\sb{1}\sp{2k+1}$ in the following way:
\begin{multline*}
\{  a, b \} \quad \cup  \quad \{  2b-a \} \quad  \cup  \quad
\{ 3b-2a, 4b-3a, \dots , kb-(k-1)a, (k+1)b-ka  \}  \quad \cup \\
\quad \cup \quad \{  (k+2)b-(k+1)a, (k+3)b-(k+2)a, \dots ,
(2k-1)b-(2k-2)a, 2kb-(2k-1)a  \}
\end{multline*}
The set $\{ 3b-2a, 4b-3a, \dots , kb-(k-1)a, (k+1)b-ka  \}$
contains the colors that are introduced from the performance of
type III Reidemeister moves $i=2$ through $k$; these colors were
already in the diagram so after these $k-1$ moves there are two of
each. The set $\{  (k+2)b-(k+1)a, (k+3)b-(k+2)a, \dots ,
(2k-1)b-(2k-2)a, 2kb-(2k-1)a  \}$ contains the colors that are
removed  from the performance of type III Reidemeister moves $i=2$
through $k$. Only one of each of the colors in this latter set was
in the coloring of the diagram before the performance of these
$k-1$ moves. Then, after these $k-1$ moves, there are none of
these colors in the coloring of the resulting diagram. In this
way, after the type III Reidemeister move corresponding to $i=k$
is performed there are $k-1=\# \{ (k+2)b-(k+1)a, (k+3)b-(k+2)a,
\dots , (2k-1)b-(2k-2)a, 2kb-(2k-1)a \}$ less colors than in the
original diagram i.e., there are only $2k+1-(k-1)=k+2$ colors
left. The result follows. (Note, in passing, that, as of this
$i=k$ type III Reidemeister move, the colors removed are the ones
that there are in two's and the colors introduced are the ones
that were previously removed, so as of this move the number of
colors in the colorings of the diagrams increases). $\hfill
\blacksquare $

\bigbreak

We remark that, for prime $p=2k+1$ and $r$ divisible by $p$ then
\[
\mathrm{mincol}\sb{r}T(2, p) \leq k+2
\]
since $\left\{  0, \; \frac{r}{p} , \; 2\: \frac{r}{p} ,\: \dots
\: , \; (p-1)\: \frac{r}{p} \right\}$ endowed with the $a \ast b =
2b-a$ (mod $r$) operation is algebraically closed. It is, in fact,
an algebraic closed substructure of $\{ 0, 1, \dots , r-1  \}$
endowed with the same operation. Moreover, if $(n, r)>1$ and
$\langle n, r \rangle = p\: (=2k+1$, for some positive integer
$k)$, then we regard $T(2, n)$ as the closure of a stacking of
$\sigma\sb{1}\sp{p}$'s, each one of which is non-trivially colored
using colors from $\left\{  0, \; \frac{r}{p} , \; 2\: \frac{r}{p}
,\: \dots \: , \; (p-1)\: \frac{r}{p} \right\}$. Teneva reduction
(Proposition \ref{prop : 2k+1}) can now be applied to each
$\sigma\sb{1}\sp{p}$ to reduce the number of colors from $2k+1$ to
$k+2$. In this way:

\begin{prop}\label{prop : p=2k+1} Suppose $r$ and $n$ are positive
integers such that $\langle n, r \rangle = 2k+1$, for some integer
$k>1$. Then:
\[
\mathrm{mincol}\sb{r}T(2, n) \leq k+2
\]
\end{prop} Proof: Omitted. $\hfill \blacksquare $

\begin{def.}
We call Teneva transformation a finite sequence of moves on knot
diagrams endowed with colorings as described in Lemma \ref{lem :
remint}. This transformation introduces some colors and removes
other colors in the colorings. We call this transformation Teneva
reduction $($Proposition \ref{prop : 2k+1}$)$ when the net effect
of the Teneva transformation is to decrease the number of colors
used in the coloring.
\end{def.}

Note that a Teneva transformation can be applied to a portion of a
knot diagram endowed with a coloring such that this portion of the
diagram looks like a $\sigma\sb{1}\sp{n}$ ($\sigma\sb{1}\in
B\sb{2}$). If the net effect of the Teneva transformation on the
braid-like portion of the diagram is to decrease the number of
colors used in the whole diagram, we call it also Teneva reduction
(see Section \ref{sect:applications}). The number of type III
Reidemeister moves should be adapted to each colored knot diagram
whose number of colors we want to reduce in order to maximize this
reduction. For the case of the $T(2, n)$ endowed with a
non-trivial $r$-coloring  such that $\langle n, r \rangle = 2k+1$,
$k-1$ of those type Reidemeister moves maximize the reduction.

\subsection{The $\langle n, r \rangle = 5$ instance.}\label{subsect:5}

\noindent

Our proof of this instance of Proposition \ref{prop : p=2k+1} relies
on understanding how the multiplication table with respect to the
$a\ast b := 2b-a \, \mod r$ operation (\cite{dJoyce, sMatveev}) of
a subset of three distinct elements from $\{ 0, 1, \dots , r-1 \}$
can be realized for a general integer $r>2$. For each of the
realizations we will then inquire into whether this subset with
the indicated table can give rise to a coloring of a $T(2, n)$ or
not.

Suppose we are given distinct $a, b, c \in \{0, 1, 2, \dots , r-1
\}$. We construct their multiplication table by considering the
distinct solutions of the equation $2x-y=z$ in $\{ a, b, c \}$. We
start by remarking that the equation $2x-y=x$ has only solutions
satisfying $x=y$ and will thus be systematically discarded. We
remark also that equalities and belonging to sets will be
understood modulo $r$.

\begin{enumerate}
\item\label{it:01.} Suppose first that the equation $2x-y=z$ has a
solution where no two variables take on the same value say, $2b-a
= c$, possibly after relabelling.
\begin{enumerate}
\item Assume further that $2a-c=b$ also holds. Then, adding these
two expressions we obtain $2c-b=a$. The multiplication table is
then shown in Table \ref{Ta:2b-a = c, 2c-a=b}.
\begin{table}[h!]
\begin{center}
    \begin{tabular}{| c || c | c | c |}\hline
                &   $a$     &      $b$  &     $c$  \\ \hline \hline
            $a$ &   $a$     &      $c$  &     $b$    \\ \hline
            $b$ &   $c$     &      $b$  &     $a$   \\ \hline
            $c$ &   $b$     &      $a$  &     $c$   \\ \hline
    \end{tabular}
\caption{Multiplication table for the $2b-a = c$, and $2c-a=b$
case}\label{Ta:2b-a = c, 2c-a=b}
\end{center}
\end{table}
Moreover, substituting the last expression in either of the first
two we obtain
\[
3(b-a)=0
\]
which implies that $3|r$ in order to be possible for $b$ and $a$
to be distinct. \item Assume now there is only one solution of
$2x-y=z$ (modulo permutation of $y$ and $z$) with no two variables
taking on the same value. This solution is $2b-a=c$.
\begin{enumerate}
\item Assume then $2a-c=c$. Then $2(a-c)=0$ which implies that $r$
is even and $c=a+\frac{r}{2}$. If $2a-b=b$ then again
$b=a+\frac{r}{2}=c$ which is a contradiction ($a, b, c$ are
distinct). Then $2a-b \notin \{ a, b, c  \}$. Further,
\[
2c-b=2\biggl( a+\frac{r}{2}\biggr) -b=2a-b \notin \{ a, b, c  \}
\]
The multiplication table is then shown in Table \ref{Ta:2a-b = X}.
\begin{table}[h!]
\begin{center}
    \begin{tabular}{| c || c | c | c |}\hline
                &   $a$     &      $b$  &     $c$  \\ \hline \hline
            $a$ &   $a$     &           &         \\ \hline
            $b$ &     X     &      $b$  &     X   \\ \hline
            $c$ &           &           &     $c$   \\ \hline
    \end{tabular}
\caption{Multiplication table for the $2a-b$ and  $2c-b$ not in
$\{ a, b, c \}$ case}\label{Ta:2a-b = X}
\end{center}
\end{table}
In Table \ref{Ta:2a-b = X} an ``X'' means the corresponding entry
does not belong to $\{ a, b, c \}$ and an empty entry means the
specification of that entry is not relevant. \item Finally assume
$2a-c \notin \{ a, b, c\}$. Then $2a-c \neq c$ which implies that
$2c-a \notin \{ a, b, c\}$. Further, if $2c-b=b$ then subtracting
this expression from $2c-a =e \notin \{ a, b, c\}$ implies $2b-a=e
\notin \{ a, b, c\}$ which contradicts the standing assumption
$2b-a=c$. Then $2c-b \notin \{ a, b, c\}$. Also, if $2a-b=b$ then
$2b-a=a\neq c$ so $2a-b \notin \{ a, b, c   \}$. The
multiplication table for this case looks like Table \ref{Ta:2a-b =
X}.
\end{enumerate}
\end{enumerate}
\item Assume now there is no solution of $2x-y=z$ in $\{ a, b, c
\}$ where no two variables take on the same value.
\begin{enumerate}
\item Suppose the equation $2x-y=y$ has solutions satisfying
$x\neq y$. Then $r$ is even. If it has more than one such
solutions (modulo permutation of $x$ and $y$), say
$2(b-a)=0=2(b-c)$, possibly after relabelling, then
$a=b+\frac{r}{2}=c$ which is impossible. Then there is at most one
solution of $2x-y=y$ with $x\neq y$ (modulo permutation of $x$ and
$y$). Suppose it is $2a-c=c$. Then $2a-b \notin \{ a, b, c \}$ and
$2c-b \notin \{ a, b, c  \}$. The multiplication table for this
case looks like Table \ref{Ta:2a-b = X}. \item If there are no
solutions of  $2x-y=y$ satisfying $x\neq y$ then the
multiplication table again looks like Table \ref{Ta:2a-b = X}.
\end{enumerate}
\end{enumerate}

\bigbreak

There are then two possibilities for $r$-colorings involving
exactly three distinct colors $a, b, c$. Either when two
distinct colors meet at a crossing the third one emerges from that
crossing. This possibility corresponds to Table \ref{Ta:2b-a = c,
2c-a=b} and further implies that $3|r$, as it was seen above. Or there is
always a color which cannot be assigned to an under-arc which ends
at a crossing. This possibility corresponds to Table \ref{Ta:2a-b
= X} where the special color is $b$. But with knots (or
non-splittable links) any under-arc ends at a crossing. These
knots and links can therefore only be colored with three distinct
colors complying with the former possibility. This is the case
with the $T(2, n)$ since for odd $n$ they are knots and for even
$n$ they are links made of two linked components. We can thus
state:

\begin{prop}\label{prop:3col} Assume $K$ is non-splittable.
If $3\nmid r$, there is no $r$-coloring of a diagram of $K$ which
uses exactly three colors from $\{ 0, 1, \dots , r-1 \}$.
\end{prop} Proof: Omitted. $\hfill \blacksquare $

\begin{prop}If $\langle n, r \rangle = 5$ then
\[
\mathrm{mincol}\sb{r}T(2, n) = 4
\]
\end{prop} Proof: We recall that, via the proof of Proposition
\ref{prop : mincor}, there is a coloring that uses exactly $4$
colors from $R\sb{r}$ whenever $\langle n, r \rangle = 5$.
Corollary \ref{cor:mincolcor2} shows us that
$\mathrm{mincol}\sb{r}T(2, n)$ here cannot be $2$. We will now
show that $\mathrm{mincol}\sb{r}T(2, n)$ cannot be $3$ thus
concluding the proof. Since $T(2, n)$ is either a knot or a link
with two linked components then an $r$-coloring involving exactly
three colors corresponds to the situation described by Table
\ref{Ta:2b-a = c, 2c-a=b} according to the discussion preceding
Proposition \ref{prop:3col}. In particular, $3|r$. Also, these
three colors are algebraically closed under the given operation.
This implies that if there is a diagram of $T(2, n)$ which is
colored by these three colors then any other diagram of $T(2, n)$
is colored by these three colors. On the other hand, the closure
of $\sigma\sb{1}\sp{n}$ is colored by three colors only if $3|n$ (adapt the argument used in the end of the proof of Proposition \ref{prop : mincol3}).
But if $3|r$ and $3|n$ then $\langle n, r \rangle \neq 5$, which
contradicts the hypothesis. This concludes the proof.

 $\hfill \blacksquare $

\begin{cor} Suppose $\langle n, r \rangle = 2k+1$ with $k>1$. Then
\[
3 < \mathrm{mincol}\sb{r}T(2, n)
\]
\end{cor} Proof: Omitted.
 $\hfill \blacksquare $

\section{Illustrative examples of Teneva reduction on rational
knots.}\label{sect:applications}

\noindent

We will now consider applications of Teneva reduction to rational
links. For facts and notation on rational knots we refer to
\cite{kl-t}.

\begin{figure}[h!]
    \psfrag{z}{\huge $0$}
    \psfrag{o}{\huge $1$}
    \psfrag{t}{\huge $2$}
    \psfrag{th}{\huge $3$}
    \psfrag{f}{\huge $4$}
    \psfrag{fi}{\huge $5$}
    \psfrag{si}{\huge $6$}
    \centerline{\scalebox{.50}{\includegraphics{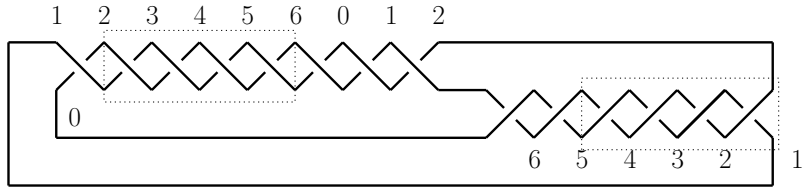}}}
    \caption{Non-trivial $7$-coloring on $N\bigl[ [8], [-6]\bigr] $
    using the $7$ colors}\label{Fi:rational86}
\end{figure}

We consider the rational knot $N\bigl[ [8], [-6]\bigr] $ endowed
with a non-trivial $7$-coloring, using all $7$ colors, view Figure
\ref{Fi:rational86}.

The two portions of Figure \ref{Fi:rational86} boxed by dotted
lines can be regarded as $\sigma\sb{1}\sp{4}$ - the one on the
left - and $\sigma\sb{1}\sp{-4}$ - the one on the right.  In
Figure \ref{Fi:tron86} we show an instance of Teneva reduction
applied to this $\sigma\sb{1}\sp{-4}$:
\begin{figure}[h!]
    \psfrag{z}{\Large $0$}
    \psfrag{o}{\Large $1$}
    \psfrag{t}{\Large $2$}
    \psfrag{th}{\Large $3$}
    \psfrag{f}{\Large $4$}
    \psfrag{fi}{\Large $5$}
    \psfrag{si}{\Large $6$}
    \centerline{\scalebox{.50}{\includegraphics{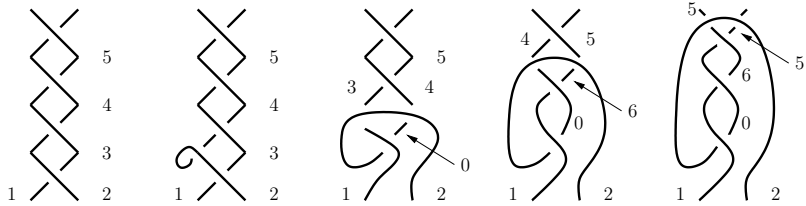}}}
    \caption{An instance of Teneva reduction}\label{Fi:tron86}
\end{figure}
Note that had we continued performing type III Reidemeister moves
in the Teneva reduction we would have increased the number of the
colors in the coloring of the rational knot.

\begin{figure}[h!]
    \psfrag{z}{\huge $0$}
    \psfrag{o}{\huge $1$}
    \psfrag{t}{\huge $2$}
    \psfrag{th}{\huge $3$}
    \psfrag{f}{\huge $4$}
    \psfrag{fi}{\huge $5$}
    \psfrag{si}{\huge $6$}
    \centerline{\scalebox{.50}{\includegraphics{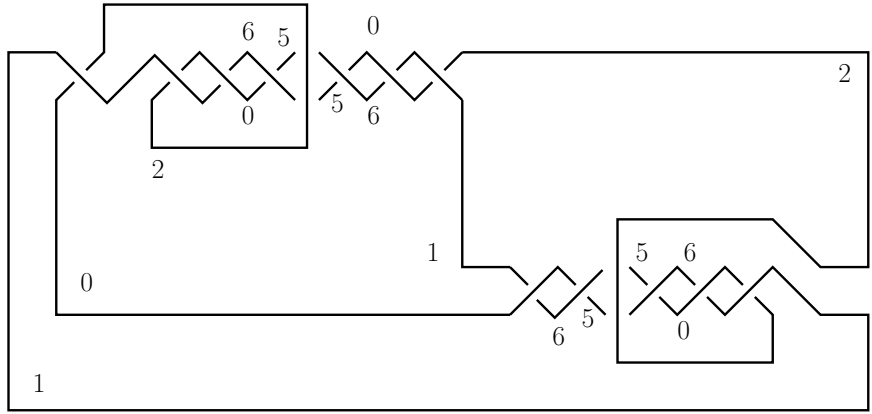}}}
    \caption{Reducing the colors of the non-trivial coloring from $7$
    to $5$}\label{Fi:trrational86}
\end{figure}

We remark that Teneva reduction can be applied in a similar way to
the $\sigma\sb{1}\sp{4}$ in Figure \ref{Fi:rational86}. With these
reductions the coloring on $N\bigl[ [8], [-6]\bigr] $ depicted in
Figure \ref{Fi:rational86} has changed to the equivalent coloring
shown in Figure \ref{Fi:trrational86}, using only $5$ colors -
less $2$ than the $7$ colors used in Figure \ref{Fi:rational86}.

\begin{figure}[h!]
    \psfrag{z}{\huge $0$}
    \psfrag{o}{\huge $1$}
    \psfrag{t}{\huge $2$}
    \psfrag{th}{\huge $3$}
    \psfrag{f}{\huge $4$}
    \psfrag{fi}{\huge $5$}
    \psfrag{s}{\huge $6$}
    \psfrag{se}{\huge $0$}
    \psfrag{e}{\huge $8$}
    \psfrag{n}{\huge $9$}
    \psfrag{ten}{\huge $10$}
    \psfrag{nt}{\huge $19$}
    \psfrag{te}{\huge $28$}
    \psfrag{ths}{\huge $37$}
    \psfrag{fs}{\huge $46$}
    \psfrag{ff}{\huge $55$}
    \psfrag{sf}{\huge $64$}
    \centerline{\scalebox{.50}{\includegraphics{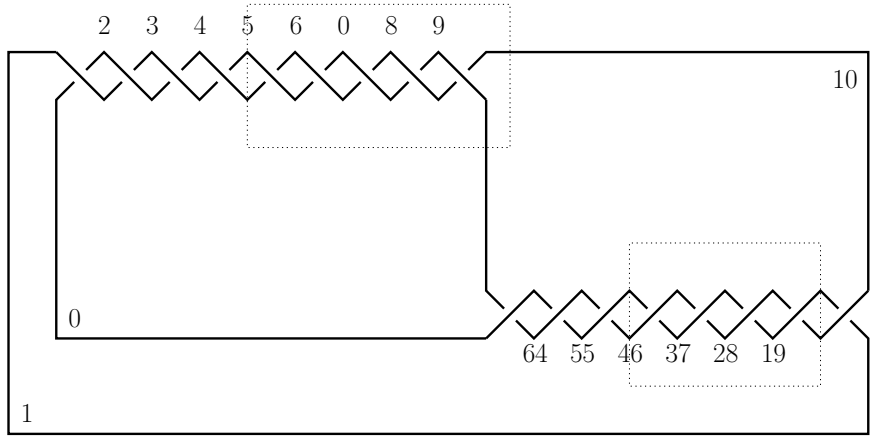}}}
    \caption{Non-trivial $73$-coloring on $N\bigl[ [8], [-9]\bigr] $
using the $17=8+9$ colors}\label{Fi:rational98}
\end{figure}

\begin{figure}[h!]
    \psfrag{z}{\huge $0$}
    \psfrag{o}{\huge $1$}
    \psfrag{t}{\huge $2$}
    \psfrag{th}{\huge $3$}
    \psfrag{f}{\huge $4$}
    \psfrag{fi}{\huge $5$}
    \psfrag{s}{\huge $6$}
    \psfrag{se}{\huge $0$}
    \psfrag{e}{\huge $8$}
    \psfrag{n}{\huge $9$}
    \psfrag{ten}{\huge $10$}
    \psfrag{nt}{\huge $19$}
    \psfrag{te}{\huge $28$}
    \psfrag{ths}{\huge $37$}
    \psfrag{fs}{\huge $46$}
    \psfrag{ff}{\huge $55$}
    \psfrag{sf}{\huge $64$}
    \centerline{\scalebox{.50}{\includegraphics{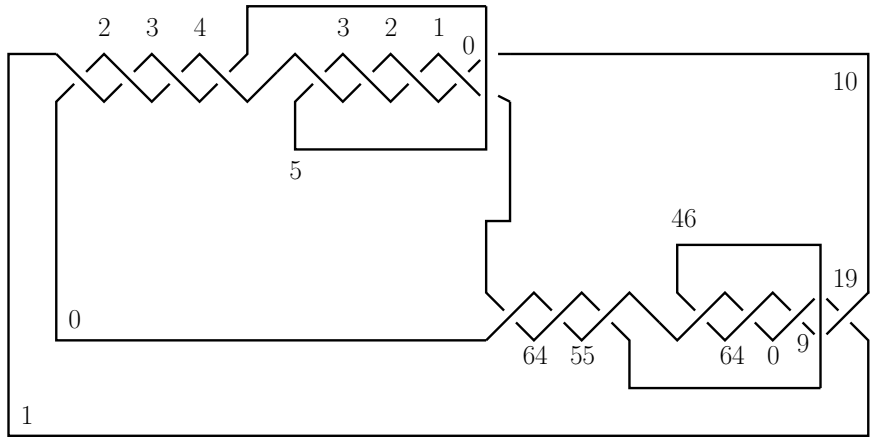}}}
    \caption{Non-trivial $73$-coloring on $N\bigl[ [8], [-9]\bigr] $
    using $12$ colors}\label{Fi:ttrational98}
\end{figure}

We now show an example of a rational knot, $R$, of prime determinant, $p$,
such that its minimum number of colors modulo $p$
 and over all diagrams is strictly less than the minimum number of
colors over minimal diagrams. This rational knot is $N\bigl[ [8],
[-9]\bigr] $. Its determinant is the prime $73(=9\cdot 8 +1)$. The
Teneva reduction in this case is similar to the preceding one so
we will just show the initial diagram endowed with a non-trivial
$73$-coloring, the portions of the diagram that will undergo
Teneva reduction (inside the boxed areas), see Figure
\ref{Fi:rational98}; and the final diagram after Teneva reduction
has been performed on the indicated portions, see Figure
\ref{Fi:ttrational98}.

In particular, note that the minimal alternating diagram on Figure
\ref{Fi:rational98} uses $17$ colors to produce a non-trivial
$73$-coloring. This is in accordance with the Kauffman-Harary
Conjecture (Conjecture \ref{conj:KH}) since this diagram has
$9+8=17$ arcs. Notwithstanding, the diagram in Figure
\ref{Fi:ttrational98} is Reidemeister equivalent to the preceding
one and uses only $12$ colors to produce a non-trivial
$73$-coloring. In this way, our Conjecture (see Section
\ref{sect:intro}) does not apply for general rational knots with
``$n$'' replaced by ``Determinant of $R$''. Moreover, this example
shows that, for rational knots, the minimum number of colors
predicted by the Kauffman-Harary conjecture, which considers only
minimal diagrams, can be further decreased using non-minimal
diagrams.

\section{Final remarks}\label{sect:final}

\noindent

The main results of this paper are the expression of the number of
$r$-colorings of a torus knot $T(2, n)$ for any integer $r>2$ and
the reduction of the upper bounds of the minimum number of colors
necessary to produce a non-trivial $r$-coloring using Teneva
transformations. A first question we would like to answer concerns
the truth of our conjecture on the minimum number of colors for
the $T(2, n)$'s.

We regard this work as a case study which motivates us to attack
more general situations. In this way, we would like to consider
other classes of knots, other classes of labelling quandles,
inquiring into the possibility of other sorts of transformations
that allow one to reduce the number of colors of a non-trivial
coloring, finding techniques to compute the minima of these colors
given a diagram of a knot (at least for some classes of knots),
etc. We aim to address these topics in future work.

\subsection{Acknowledgements} \label{subsect:ack}

\noindent

It gives L.K. pleasure  to thank the National Science Foundation
for support of this research under NSF Grant DMS-0245588.

 P.L. acknowledges support by {\em Programa Operacional
``Ci\^{e}ncia, Tecnologia, Inova\c{c}\~{a}o''} (POCTI) of the {\em
Funda\c{c}\~{a}o para a Ci\^{e}ncia e a Tecnologia} (FCT)
cofinanced by the European Community fund FEDER. He also thanks
the staff at IMPA and especially his host, Marcelo Viana, for
hospitality during his stay at this Institution. \bigbreak


\begin{thebibliography}{99}


\bibitem{aps}
    M. Asaeda, J. Przytycki, A. Sikora,
    Kauffman-Harary conjecture holds for Montesinos knots,
    J. Knot Theory Ramifications 13 (4) (2004) 467-477



\bibitem{Birman}
        J. Birman, Braids, links, and mapping class groups,
        Annals of Math. Studies 82, Princeton University Press,
        Princeton, NJ, 1974





\bibitem{DL}
        F. M. Dion\'\i sio, P. Lopes,
        Quandles at finite temperatures II,
        J. Knot Theory Ramifications, 12 (8) (2003) 1041-1092




\bibitem{Fox}
    R. H. Fox, A Quick Trip Through Knot Theory,
    in M. K. Fort, Jr. (Ed.), Topology of 3-Manifolds and Related Topics, Georgia, 1961,
    Prentice-Hall, 1962, 120-167




\bibitem{hk}
        F. Harary, L. H. Kauffman,
        Knots and graphs I. Arc graphs and colorings,
        Adv. in Appl. Math. 22 (3) 312-337 (1999)





\bibitem{dJoyce}
        D. Joyce,
        A classifying invariant of knots, the knot quandle,
        J. Pure Appl. Alg. 23  37-65 (1982)








\bibitem{lhKauffman}
        L. H. Kauffman, Knots and physics, third ed.,
        Series on Knots and Everything 1,
        World Scientific Publishing Co., River Edge, NJ 2001



\bibitem{kl-t}
        L. H. Kauffman, S. Lambropoulou,
        On the classification of rational tangles,
        Adv. in Appl. Math. 33 (2) 199-237 (2004)






\bibitem{pLopes}
        P. Lopes, Quandles at finite temperatures I,
        J. Knot Theory Ramifications 12 (2) 159-186 (2003)

\bibitem{sMatveev}
        S. V. Matveev, Distributive groupoids in knot theory,
        Math. USSR Sbornik 47 (1) 73-83 (1984)







\end{thebibliography}
\end{document}